# 三角剖分图染色与换色通道

甘润东　　南开大学　　ynztgrd@163.com

# 摘要

三角剖分图染色，是平面图染色的充分情况。本文将重点研究三角剖分图，以及染色工具换色通道的性质。通过构造，将原来顶点的 4 染色，变换为边的 3 染色，和三角形的 2 染色，从而将等价的染色方案合并。换色通道利用边染色的可传递性，体现出三角剖分图的整体性。根据换色通道的性质得到两条约束条件，让换色有规律可以遵循。最终使得在处理约化 5 次顶结构上获得突破，证明了该结构的可约性，从而否定了五色图最小反例的存在。

**关键词**：三角剖分图；肯普换色链；换色通道；不友好情况

**MR(2000)分类号**：05 C 15（平面图染色）

# 目录

# 第一章、定理的阐述与历史

## 1.1 四色定理的阐述

四色定理，又称为四色猜想，四色问题。（可称为 4CT：4 color theory）。

4CT 是指，任意的一个平面图，都可以用 4 种颜色区分开来，使得相邻的两个区域颜色不同。

进一步的，我们可以将平面图作对偶化处理，添加更多的连接边，使得所有顶点尽可能的相互连接。这样使得其成为一个三角剖分图。那么显然的，如果任意的三角剖分图都是四色的，则它的任意子图也是四色的。

因此，我们的研究，就以三角剖分图作为对象，不再从其他图作为研究的切入点。

## 1.2 四色定理的发展历史

### 1.2.1 早期发展

4CT 于 1852 年由 Francis Guthrie 提出。

1878 年，Alfred Kempe 提出肯普链，宣布证明成功。

11 年后，希伍德 P.J.Heawood 提出反例，指出肯普证明的漏洞。

1880 年，Peter Guthrie Tait 以边染色为切入点，借助哈密顿路径，宣称证明 4CT。

1946 年，Tuttle 提出反例，Tait 的证明不成立。

尽管肯普和希伍德均未能破解四色问题,但是,他们二人均为四色问题的最终获证乃至图论的发展做出了早期贡献。尤其是肯普所带来的不可约集和不可避免集的概念。

### 1.2.2　计算机证明

1976 年,4CT 由阿佩尔 K.Appel 与哈肯 W.Haken 以机器方法证明。

1944 年，西缪尔 P.Seymour 等人优化了证明。缩短了时间和证明篇幅。

2005 年，贡蒂埃 Georges Gonthier 给出形式证明（formal proof）。

## 1.3　研究意义

时至今日，尚未出现四色定理的人工证明。某种意义上，四色定理体现的是二维平面的一个固有属性，但更深入的，却是和 NP 完全问题有一定的相关性。解决四色定理的人工证明，有助于理解平面的性质，有助于处理 NP 问题。这个过程中产生的新工具，有助于图论中染色理论的发展。

# 第二章、前期准备与新工具

## 2.1 准备工作

在证明工作开始之前，我们先做一些前期的准备工作。解决和定义顶点染色与边染色。

**2.1.1 定义** 在三角剖分图的染色中，顶点染色用"1、2、3、4、5"来指代五种颜色，边染色用"a、b、c"来指代三种颜色，三角形面染色用"上、下"来指代两种颜色

**2.1.2 引理** 若顶点数为 n 的三角剖分图 G(n) 是 4 色的，则其中每个三角形都是由 a、b、c 三色边构成。

证明：因为 G(n) 是 4 色的，则设顶点染色的四种颜色分别是 1、2、3、4。将其置于色位图中。

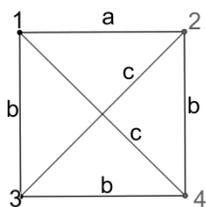

图 **2.1.2**：色位图，将顶点染色放在固定的位置

设定横向边染色为 a，纵向边染色为 b，对角边染色为 c。显然的，三角剖分图 G(n) 中，任意三角形的顶点染色，都是 1234 颜色中取 3 个。通过图 1.能够直观得出结论，任意三角形必然有 a、b、c 三种颜色的边。

**2.1.3 推论** 由于每个三角形都是 abc 三色边，然而 abc 的顺

序存在顺逆时针的区别。所以顺时针的 abc 三角形为"↑"，逆时针的 acb 的三角形为"↓"。

2.1.4　**定义**　在一个已经成功染色的三角剖分图 G(n)中，我们把一条迹 J(m)（迹是无重复边的通路）中的边染色作计数。那么有 J(m)=A+B+C。色边的边数如果是奇数，记为 A1 或 B1 或 C1，如果是偶数记为 A2 或 B2 或 C2。如果迹中没有出现这种颜色的边，算作是偶数边。

故而相同颜色边存在计算关系：

A1+A1=A2　　　A1+A2=A1　　　　A2+A2=A2

不同颜色边数不存在计算关系。

例如 J(m)=A1+B2+C1+A1=A2+B2+C1

2.1.5　**定理**　在一个已经成功染色的三角剖分图 G(n)中，任意一条闭合的迹，都满足 a、b、c 三条边的数量奇偶性一致。即 A、B、C 同奇或者同偶。

证明：因为三角形都是 abc 三边，所以存在路径替换的关系：

a=b+c　　b=a+c　　c=a+b　　b+c=a　　a+c=b　　a+b=c

路径替换不改变染色方案。

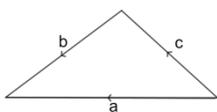

图 2.1.5：路径替换

任意一条闭合的迹,内部都将包含若干个三角形，通过路径替换，不断的减少内部包含的三角形，直到只剩下一

个三角形 abc。此时 J(abc)=A1+B1+C1。根据其逆过程，可以显而易见的得出 ABC 奇偶性一致的结论。

当我们做到这一步时，逐渐的让三角剖分图的染色，变得可以计算了。

2.1.6　引理　在一个已经成功染色的三角剖分图 G(n)中，任意一个轮图 L(m)，都满足 L(m)=0。L(m)等于各个三角形面染色的相加。相加满足以下计算法则。

↑＋↓＝0　　　　↑＋↑＋↑＝3↑＝0　　　　↓＋↓＋↓＝3↓＝0

且满足交换律和结合律。

例 L(m)=↑＋↑＋↓＋↑＋↑＋↓＋↑＝5↑＋2↓＝3↑＝0

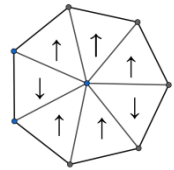

图 2.1.6：轮图中的面染色关系

2.1.7　证明　设中心顶点的颜色为 1，轮上颜色为 234。利用色位图辅助。

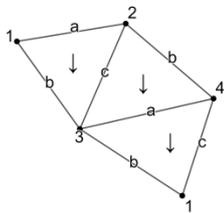

图 2.1.7：轮图中的面染色关系

显然的，轮图上的颜色符合在三角形 234 上的顺逆时针变化。因为 n 边形的轮图，经过 n 个 ↑↓三角形的变化后，轮上端点颜色是一致的。所以可以明显看出，顺时针或者逆时针三次变化后，顶点回到最初的颜色。也就是 3↑＝0。顺一次逆一次，顶点颜色还原。也就是 ↑＋↓＝0。

事实上，这些结论都非常的显而易见，以至于不需要单独的证明一遍。

2.1.8　总结　不论是边染色还是面染色，都没有改变什么。只是让我们能够更清晰的了解问题，这类似于代数中的换元思想。通过换元简化问题的难度。

## 2.2　换色通道

工欲善其事，必先利其器。一个高效、准确、清晰、明显的工具，可以极大的简化图论中的工作。使得问题由局部性拓展到全局性，避免在局部有解但全局上出现矛盾。

如果一个图 G(n)是 4 色的，且至少有一种染色方案时，就意味着另外等价的 11 种染色方案。共计 12 种染色方案。这取决于四种颜色间的轮换关系。这个些多余但是等价的染色方案，极大的增加了工作的复杂程度。于是把它处理约化为更简单的形式。

2.2.1　推论　图 G(n)是一个顶点可以 4 色的三角剖分图。

存在 4-顶点染色（1、2、3、4），12 种等价染色方案。

存在 3-边染色（a、b、c），6 种等价染色方案。

存在 2-三角形面染色（↑、↓），2 种等价染色方案。

其中面染色是最根源性的表达，避免了复杂的等价染色方案。

2.2.2　**定义**　图 G(n)是一个顶点可以 4 色的三角剖分图，如果我们将 G(n)中的两条色边相互替换，会得到一个新的染色方案，与之前的染色方案等价。我们将这一过程定义为一次边变换，例如 a-c 变换，记做 a-c，a-c 变换等价于 c-a 变换。同样的，存在 a-b 变换，b-c 变换。

2.2.3　**定义**　将一次边变换的过程动态处理，我们会得到至少一条首尾相连的变换路径。我们称这条由三角形构成的路径为换色通道。例如 a-c 边变换形成的通道，就称为 a-c 换色通道。

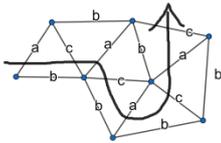

图 2.2.3：a-c 换色通道

证明：因为所有的三角形都是 abc 三边，且每条边都是两个三角形的公共边。所以一个三角形上的边颜色互换，就会依次影响到下一个三角形，重复这个过程，直到首尾相连时停止。此时这条换色路径就是换色通道。显然的，通道内部是两种颜色边交替出现，边界是第三种颜色。例如 a-c 通道中，a、c 交替出现，b 作为通道边界。

2.2.4　**注**　边变换可能有多条换色通道，因为可能存在纯色闭合边，使得另外同类的换色通道被隔离开来。

2.2.5　**推论**　显然的，一次边变换会使得变换路径上的所有三角形 ↑↓ 发生改变。

**2.2.6 引理** G(n)中一共有三种换色通道，每种通道至少一条。

**2.2.7 引理** 相同的换色通道不会交叉

证明：如果交叉，将存在一个三角形有多个同类通道通过，不能满足三角形换色路径唯一的约束条件。

**2.2.8 引理** 一条换色通道换色之后，不会造成整体染色方案的失效。

**2.2.9 重要例子** 以两个相连三角形，得出换色时的一个重要规律。这是我们以后证明的核心方法之一。

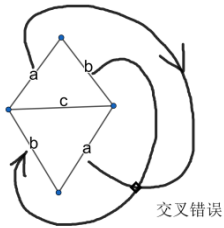

图 2.2.9：不可能出现的换色通道交叉情况

**2.2.10 推论** 以某一条边为出发，向某一个方向开始换色通道换色，其终点一定是这条边的另一个方向。有始必然有终。

## 2.3 换色通道处理两个重要反例

在四色问题历史上，出现过两个重要的反例，分别击碎了肯普和泰勒的证明方法。然而由于换色通道是更加贴近本质的方法，在这两

个反例上，体现出了很高的处理能力。在这一个章节，将使用大量的图，用以尽量阐述清楚换色通道的动态过程。

## 2.4 泰勒的哈密顿路径与塔特的反例

### 2.4.1 塔特反例图

### 2.4.2 塔特反例图的对偶图并边染色 图8.其中虚线边为a边，实线边用 b、c 交替染色。

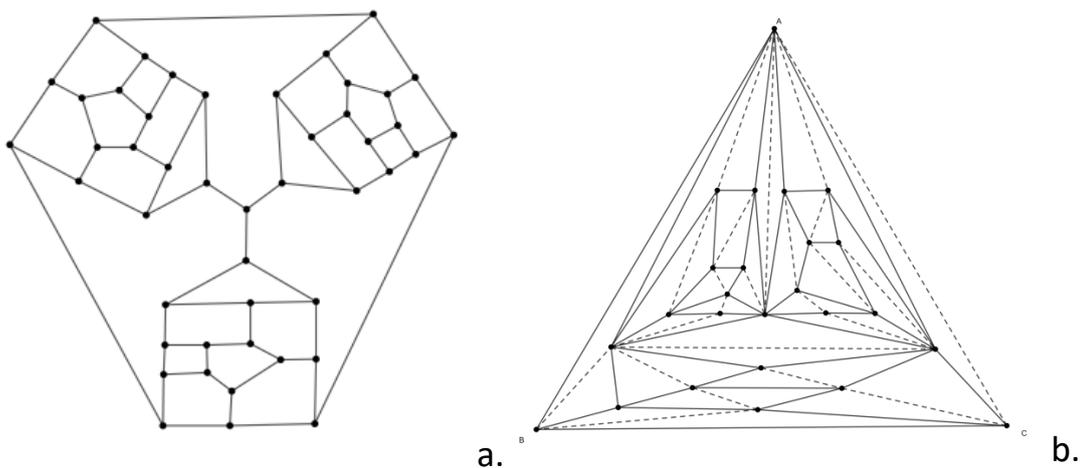

图 2.4.2:a 是塔特反例图，b 是其对偶图

### 2.4.3 证明 显然的，染色成立。每个三角形都是 abc 边。并且出现了两条换色通道，其中一条闭合的 a 边在左下角。某种意义上，换色通道类似于泰勒提出的哈密顿路径。而且由于不要求换色通道满足哈密顿路径，所以换色通道在出现闭环的问题上没有失效。

## 2.5 肯普的换色链与伍德的反例

### 2.5.1 伍德的反例图

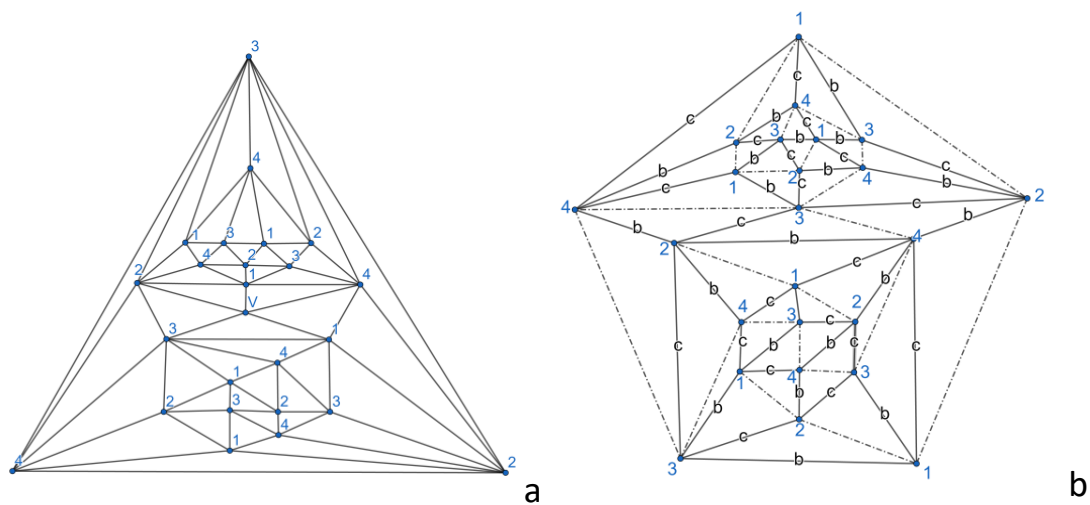

图 2.5.1：a 是伍德反例图，b 是其同构图染色

### 2.5.2 伍德的反例图进行边染色 并将内五边形外翻，形成同构图。图 b。其中虚线为 a 边。当存在 4 种颜色的时候，为 aabac 形式。

### 2.5.3 第一次换色

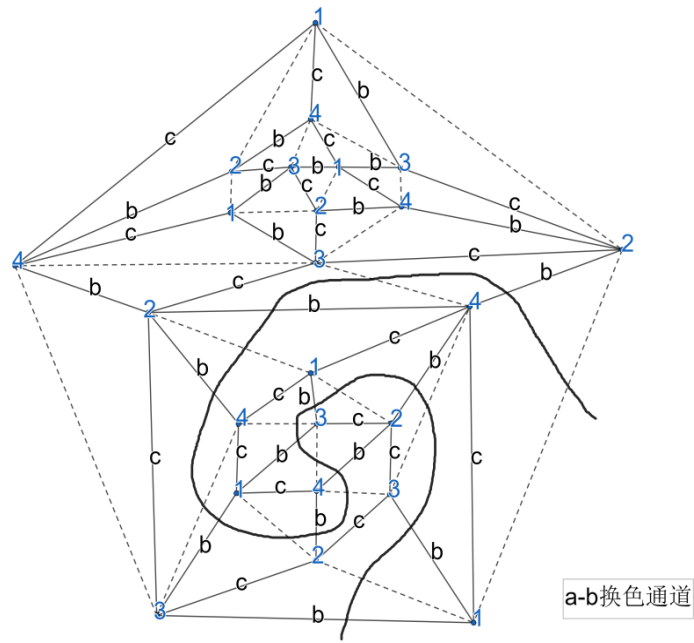

图 2.5.3：第一次换色中的一条 a-b 通道

## 2.5.4 第二次换色

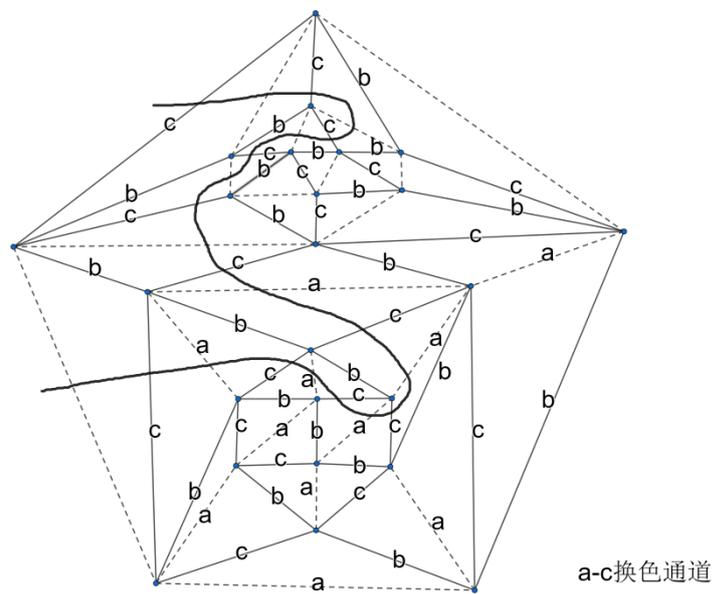

图 2.5.4：第二次换色中的一条 a-c 通道

2.5.5 **第三次换色** 换色后染色成功。得到 aaabc 的边染色情形。此时任意定义一个点的颜色，就可得到一套正确的染色方案。

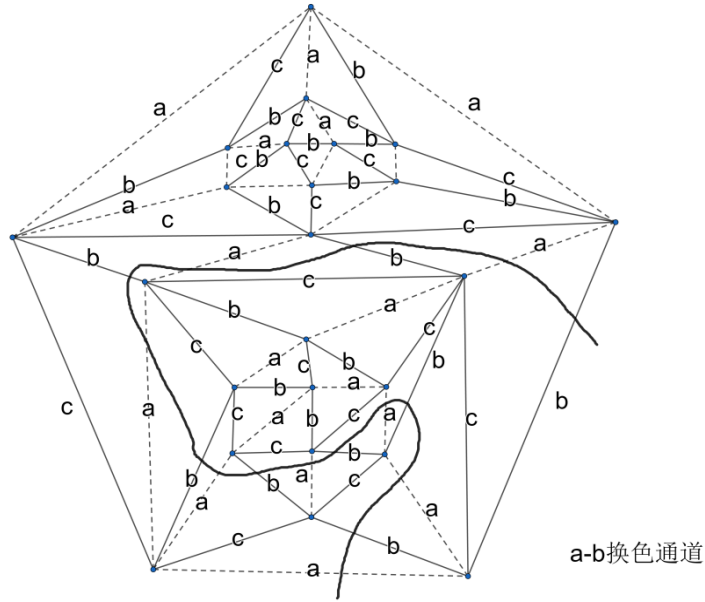

图 2.5.5：第三次换色中的一条 a-b 通道

## 2.6 分析

因为五边形的顶点染色，如果只能出现三种颜色，那么此时边染色的情况一定是 aaabc。也就是说，必然要使得三个最多的边色靠在一起。所以换色的时候，我们的目的就是尽可能的将最多的色边挪到一起去。（暂且定义最多的色边为 a 边，包括后面的章节都沿用此设定）。伍德的反例图之所以困难，是因为 a-c 换色之后，形成了新的 a-b 通道。如果使用换色链，这个情况就变得异常难处理。而换色通道就能非常直观的判断，哪些换色是可以存在，哪些是不可能存在的。从而我们可以清晰的实现换色通道的连续换色，直到出现 aaa 三连的情况，将其约化。

# 第三章、特殊情况与证明

## 3.1 肯普的困境

肯普的思路是很缜密且大部分情况是有效的，但他在遇到反例的时候，最大的困境在于，换色链太多，没有办法将换色这一个过程规范。五边形一共 5 个顶点，每个点都可以产生 3 种换色链，最多一共 15 种。在第一次尝试完这 15 种色链之后，如果换色没有成功，第二次换色将可能是 210 种方案。第三次可能接近 $15^3$ 这个大小。以至于根本无法处理连续换色，更没法证明一定能够换色成功。

显然的，造成以上困境的主要原因，是因为这么多种方案中，包含了大量重复等价的染色方案，增加了判断的难度。但在使用换色通道之后，极大的合并了等价的染色方案，并且换色的规律变得有迹可循。换色通道换色的核心约束：一是引理，换色通道不交叉。二是定理，ABC 同奇或同偶。

## 3.2 对肯普工作的复盘

3.2.1 在三角剖分图 G(n)中，有 P 个顶，q 条边，r 个面。那么 k 次顶的数目是 P（k）个，P（0）=P（1）=0。根据 Euler 公式，P-q+r=2，可知 q=6P-12。也就是说，在顶足够多的三角剖分图 G(n)中，有至多 5 次的顶存在。

若四色猜想不成立，我们可以取得一个反例，它具有最少的顶数。如果他不是三角剖分，我们可以增加边，使其成为三角剖分图 T(n)，

而其顶数不变。任意小于顶点小于 n 的平面图皆四色，而 T(n)不是。

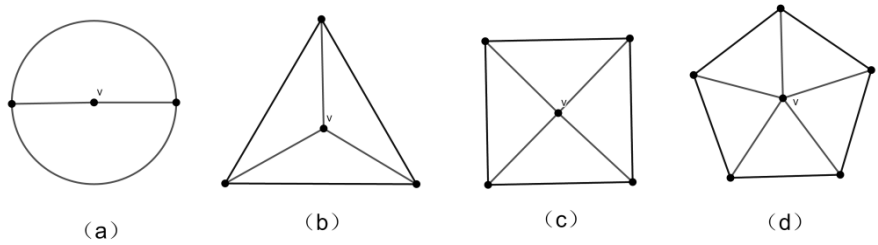

图 3.2.1：T(n)中的四种基础结构

在（a）（b）（c）（d）四种结构中，如果含有（a）（b），显然的，可以把 v 点染上相异的颜色。故而 T(n)不含有结构（a）（b）。

若 T(n)中含有结构（c），以相同的方法证明。但当顶点颜色分别是 1234 的时候，不能立刻给 v 染色。此时，肯普使用了换色链的方式，证明了（c）结构的可约。证明过程虽不复杂但不在此赘述。

在证明（d）结构的可约性上，肯普遇到了挫折，并没有证明成功，原因在第二章已经介绍过了。

### 3.2.2 换色通道处理（c）结构

如果（c）的顶点分别是颜色 1234，那么我们把边染色后得到图 a.

显然的，我们做一次关于左边的 a-b 换色通道换色。那么 a-b 通道从左边出发，由于换色通道不能交叉的约束性，必然从上方或者下方的边回归，如图 b。无论哪种情况，都形成了边为 aabb 的四边形。将顶点染色，就是 1214，于是 v 点可以染上第四种颜色。证明了（c）结构可约。

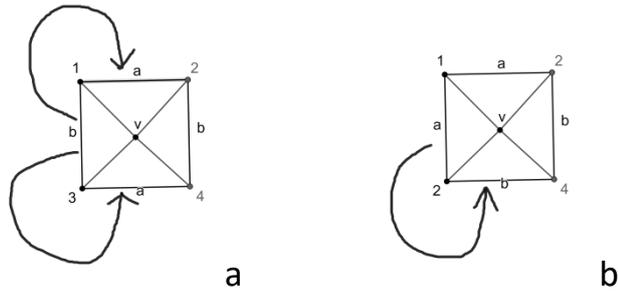

图 3.2.2：（c）结构中的换色通道和换色结果

### 3.2.3 换色通道处理（d）结构

同样的，将其边染色。如果要使得五边形上的顶点染色只用三种颜色，则必然出现 aaabc 的形式（1）。由于边染色的轮换关系，我们规定最多色边为 a 边，后面都沿用此设定。

当出现四种颜色的时候，必然是 aabac 的形式（2）。

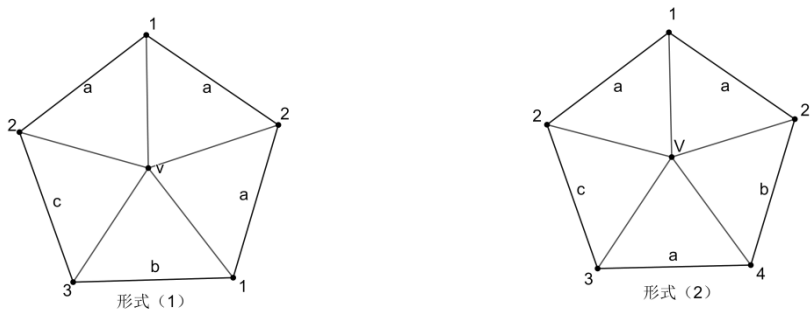

图 3.2.3：（d）结构中的两种染色形式

我们所做的，就是把（2）转变为（1）.事实上，在上一个章节，我们已经成功的实现了对于反例的突破，得到了 aaabc。现在的问题是，如何证明，一定能够通过换色通道，实现从（2）到（1）的转换。注意在这个过程中，保证永远有 3 条 a 边。

### 3.3 友好情况与不友好情况

### 3.3.1 友好情况

我们称通过换色通道，一次换色就能形成（1）形式的情况为友好情况。在肯普的换色链方法上，大部分对应的就是友好情况。通道方向用箭头表示。显然的，通道的边界将是连接两点的一条纯色边。a-c 或者 a-b 直接达成（1）形式。如图 a.

### 3.3.2 次友好情况

需要通过两次换色才能实现（1）形式。a-b 换色后再 a-c 换色。例如图 b.。

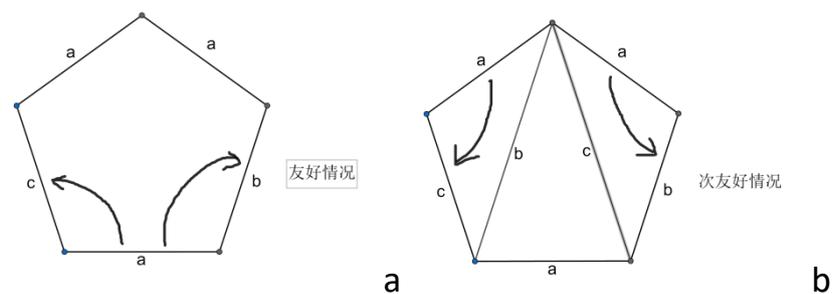

图 3.3.1：a 是友好情况，b 是次友好情况

### 3.3.3 不友好情况

a-b 或者 a-c 换色后，导致原来的换色通道受损改变，形成新的通道，以至于两次换色之后依然是（2）形式。

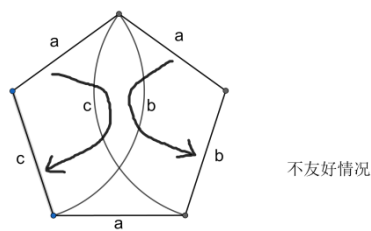

图 3.3.3：不友好情况

显然的，以上三种形式包含了所有的换色策略。a-b 和 a-c 换色不可能

有更多的选择，由于换色通道引理的限制。所以，问题的集中矛盾，放在了不友好情况上。

## 3.4 极不友好情况

不论是友好还是次友好情况，都可以达成（1）形式。而当不友好情况换色后，形成了新的（2）形式后，我们避免还原成原方案，继续向前换色。此时又将可能出现三种情况：友好，次友好，不友好。以此类推。

3.4.1 **定义** 我们称，每次换色后都是不友好的情况为极不友好情况。显然的，如果某一次换色后不再出现不友好情况，则极不友好情况就破缺了，必然出现（1）形式，实现（d）结构的约化。

在极不友好情况中，每次的不友好情况虽然都是（2）形式，但是是有区别的。同样是 aabac，但是位置不同。我们根据 a 的位置，定义指针方向。如图 3.4.1。

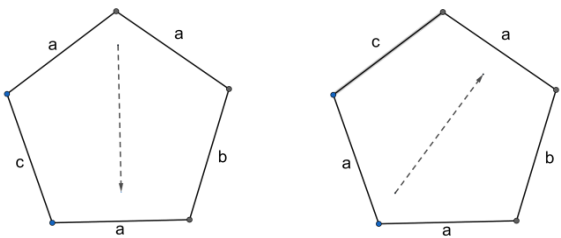

图 3.4.1：两个方向的指针

显然的，一共有 5 个方向。连续的不友好情况，将会改变指针方向。

3.4.2　推论　因为平面图是有限图，三角剖分图是有限图，所以染色方案的个数是有限的。因为每次换色都形成了区别于换色前的染色方案，所以如果极不友好情况出现，则意味着不友好染色方案必然是循环的。反之，如果证明了不循环，则必然出现循环上的破缺，也就是出现（2）转变为（1）形式。

## 3.5　极不友好情况与旋钮

### 3.5.1　构造极不友好情况

如果存在极不友好情况，那它一定是循环的。

我们不妨尝试构造一个极不友好的结构。

图 3.5.1：一个简单而对称的不友好情况

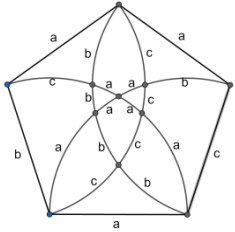

3.5.2　注　如果我们一个方向换色次数太多，会使得构造的结构非常复杂，面临巨大量的分类讨论。故而采用左右两个方向各 3 次换色，使获得指针方向统一的（2）形式。如果染色方案此时一致，则证明可以循环，也就是至少存在一种极不友好情况。如果染色方案不一致，则证明在第一次指针重合的时候不能循环。但不排除在第二次，第三次乃至第 n 次的时候，会循环。

3.5.3　结论　显然的，通过换色，让指针重合的时候，染色方案并不一致。（此处增加一次 b-c 的全局换色，协调 b 边 c 边的相对位置）。故而证明在第一次指针重合的时候不能循环。

### 3.5.4　闭合的纯色边

在图 3.5.1 中，出现了一个特殊的结构。在底边出现一个拥有 a 边围成的一个闭合区域。

显然的，闭合的纯色边，意味着两种颜色交替的色链。同时意味着，纯色边内的另外两种色边可以相互替换而不影响整体的四色性。

所以，我们将 a 圈内的 b-c 换色。使得原来的通道被改变，将不友好情况变成了友好情况。在换色通道理论下，我们称由纯色边围成的区域为旋钮。因为它就像旋钮阀门一样，改变了通道的方向。

## 3.6　旋钮

3.6.1　**定义**　一般的，我们称由纯色边围成的最小区域为旋钮。根据边的颜色，分为旋钮 a，旋钮 b，旋钮 c。旋钮之间相互独立，但可能因为其他旋钮的转动而使得自身消失。

3.6.2　**引理**　旋钮的转动，一定会改变通道属性。a-c 通道变成 a-b 通道，a-b 通道变成 a-c 通道。

当 a-b，a-c 出口不一致时，旋钮的转动会使得出口的位置互换。

3.6.3　**注**　其实我们在伍德反例图里面已将发现了旋钮，转动旋钮，马上就可以使得结构换色成功。得到 aaabc 形式。

### 3.7 边染色的奇偶性判断

#### 3.7.1 回顾将要使用的工具

在一个已经成功染色的三角剖分图 G(n)中，任意一条闭合的迹，都满足 a、b、c 三条边的数量奇偶性一致。即 A、B、C 同奇或者同偶。

另外还满足一套计算法则。把一条迹 J(m)中的边染色分别作计数。如果是奇数记为 A1 或 B1 或 C1，如果是偶数记为 A2 或 B2 或 C2。如果迹中没有出现这种颜色的边，算作是偶数边。那么有 J(m)=A+B+C。

故而相同颜色边存在计算关系：

A1+A1=A2        A1+A2=A1        A2+A2=A2

满足交换律和结合律。

故而依据此，我们可以将未明的结构中各个边的奇偶性推测并标注出来，且不失准确性。为了方便，统一暂定色边和色边数都用 a、b、c 在图中标识。

#### 3.7.2 奇偶性推断

我们将在后面的所有图中，增加边的奇偶性表示：a1、b1、c1、a2、b2、c2。

例如图 23.1。当确定 a 为 a1 的时候，自然判断出 b1 和 c1。

如图 23.2。当确定其他边的时候，判断出空缺边为 a1。

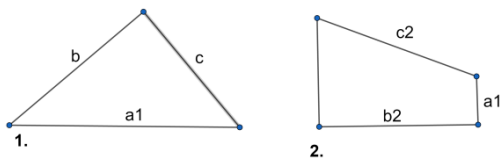

图 23.：奇偶性推断的示例

在后面的图里面，由于图

片篇幅限制，自行推断和添加带有奇偶标注的边 a1、b1、c1、a2、b2、c2。

### 3.7.3　不友好情况各边的奇偶推断

假若极不友好情况存在，则一定会包含图 3.7.3 的结构

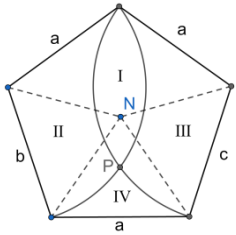

图 3.7.3：不友好情况中的 N 点在Ⅰ区域

其中有①②③④四条线，这四条线都是在某次换色时的纯色边界。其中①②线交于 P 点，③④交于 N 点。那么 N 的位置一共可以在Ⅰ、Ⅱ、Ⅲ、Ⅳ四个区域，我们对这四个区域进行分情况讨论。

### 3.7.4　注　如果 N 点在线上，则认为 N 点在两个区域，符合两个区域的奇偶性限定。实际上，N 在线上反而让问题更简单，就不再更多的做说明。

### 3.7.5　N 点在Ⅰ区域时候的情况

此时各色边的奇偶性合法，可以成立，并且我们发现存在 a 的闭合圈，也就是一个旋钮。转动旋钮，使得原有的换色通道破坏，转变成了友好情况。这在前面的章节已经得到了说明。（见图 3.5.1）

### 3.7.6　N 点在Ⅱ、Ⅲ区域时候的情况

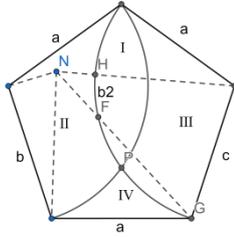

图 3.7.6：不友好情况中的 N 点在Ⅱ、Ⅲ区域

此时各色边的奇偶性合法，可以成立，并且我们暂时没有发现旋钮。但值得注意的是，这里单独出现了一个 b2 边 HF（在Ⅲ区域时时 c2 边）。因为换色通道的性质，所以出去的 b-c 通道，要么从另一边回来，从而在 HF 上得到一条 a 的纯色边界，故而依然存在旋钮。转动旋钮，使得极不友好情况破缺。要么不回来，则形成一条 HF 上出发，到右边顶点的 a 的纯色边，经过一次换色后，必然形成旋钮，使得不友好情况破缺。思路简单，但是分类讨论情况较多，篇幅所限，不一一展示。

### 3.7.7  N 点在Ⅳ区域时候的情况

在这个情况下，出现了明显的奇偶性矛盾。显然的，③④线必须在 P 点的上方通过，否则就无法满足奇偶性约束。换言之，这种情况是不能存在的。

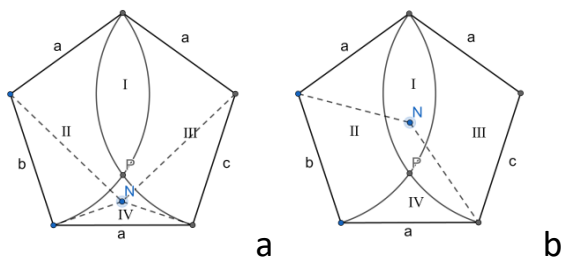

图 3.7.7：a 是不友好情况中的 N 点在Ⅳ区域。b 是 N 点的位置限定

### 3.7.8  综上四种情况所述  除了 N 点在Ⅳ区域时候不能成立，其余情况都能成立且都存在能够破坏原有换色通道的旋钮。通过转动

旋钮，使得极不友好情况破缺，得到友好情况，从而实现 aaabc 的形式，也就是能够将肯普提出的（d）结构约化，从而填上了四色证明的最后一块积木。

## 3.8  四色定理证明

至此我们讨论完了所有的情形。

若四色猜想不成立，我们可以取得一个反例，它具有最少的顶数。如果他不是三角剖分，我们可以增加边，使其成为三角剖分图 T(n)，而其顶数不变。任意小于顶数小于 n 的平面图皆四色，而 T(n)不是。在（a）（b）（c）（d）四种结构中。（a）（b）显然能够约化。（c）通过换色通道一次换色后约化。

（d）结构分情况讨论，利用换色通道：

①友好情况：一次换色后约化

②次友好情况：两次换色后约化

③不友好情况：若干次换色后约化

④极不友好情况：N 点必然在Ⅰ、Ⅱ、Ⅲ区域，从而直接或者间接的形成旋钮，转动后使得极不友好情况破缺，从而将其约化。

综上所述，反例 T(n)中不存在（a）（b）（c）（d）四种结构，这和三角剖分图的性质矛盾，故反例 T(n)不存在。

从而四色定理成立。

# 第四章、总结

四色猜想自提出以来，已经过了一百多年。在 1976 年的时候，首次通过计算机得到证明。四色问题的非计算机证明到目前为止主要证明思路，依然停留在 Tait 和 Kempe 的阶段。这两种证明思路的设计都是巧妙和慎密的，对四色问题的理论证明都做出了历史性贡献。但是在他们遇到反例时却表现得软弱无力，暴露出各自的不足。

四色问题的难度主要体现在，面对组合爆炸，传统的以集合论为思维模式的数学方法的失效。要么是将结构剥离开来，亦或是寻找图与代数之间的关系。我曾花费大量的时间，实现了三角剖分图的代数表达，并尝试用代数方法，求证其必然有解。但结果是失望的，形式上的优美，并无助于解决四色问题。

换色通道无疑是目前解决四色问题最佳的工具。它巧妙的整合了等价的染色方案，同时将整体性体现出来，让整个结构变得联动。最关键的是，推导得到几条非常重要且实用的约束条件。一是色边数的同奇同偶，二是换色通道不交叉。这使得无从下手的结构可以被调整，减少了无意义的尝试，让换色变得具有规律性，能够被分类讨论。

某种意义上，换色通道是四色问题中更本质的规律。在它的基础上我们很容易推导出 Tait 和 Kempe 的工作。也可以认为是将两者的方法结合到了一起。不过对于作者我来说，在得出换色通道这个方法之前，并没有系统的了解过这个问题的研究历史，而是在后来的整理写作中，才了解到两位前辈的方法。换色通道的成功，也无疑证明了

两位前辈在解决四色问题上，思维方向的正确性。

展望一下，四色猜想已然成为四色定理。但对于数学家来说，图论研究的事业方兴未艾。大量问题还有待解决，例如图的顶点着色、边着色、拓扑图论中的着色问题，以及这三种量的任何组合的着色问题。现代图论之父塔特(W.T.Tutte,1917-2002)感慨道，"四色问题是冰山一角、楔之尖端、孟春一啼"。而现在，在新工具的推动下，冰山正在慢慢浮起。

# 参考文献